\documentclass[final,11pt]{article}

\usepackage{geometry}                
\usepackage{graphicx}
\usepackage{amssymb}
\usepackage{epstopdf}

\usepackage{amsmath}

\newtheorem{theorem}{Theorem}

\newtheorem{remark}{Remark}
\renewcommand{\baselinestretch}{1.1}
\def\eqsp{\noalign{\vskip 5pt}}
\def\dsp{\displaystyle}
\newcommand{\bftau}{\mbox{\boldmath $\tau$}}
\newcommand{\bfeta}{\mbox{\boldmath $\eta$}}
\newcommand{\ignore}[1]{}

\newcommand{\eql}[1]{\begin{equation}\label{#1}}
\newcommand{\eqm}{\begin{eqnarray}}
\newcommand{\enm}{\end{eqnarray}}
\newcommand{\eqml}[1]{\eql{#1}\begin{array}{rcl}}
\newcommand{\enml}{\end{array}\end{equation}}
\def\grad{\nabla}
\newcommand{\eqmno}{\begin{eqnarray*}}
\newcommand{\enmno}{\end{eqnarray*}}

\usepackage{bm}
\usepackage{graphicx}
\usepackage{framed}

\usepackage[all]{xy}

\title{On jump relations of  anisotropic elliptic interface problems}
 
\author{Baiying Dong\thanks{School of Mathematics and Statistics, Ningxia University, Yinchuan, 750021, China.}
\and Xiufang Feng\thanks{School of Mathematics and Statistics, Ningxia University, Yinchuan, 750021, China.}
\and Zhilin Li \thanks{Corresponding Author: CRSC  \&
Department of Mathematics, North Carolina State University, Raleigh,
NC 27695-8205, USA. Z. Li is partially supported by a Simon Foundations grant.}
}
\date{}                                           
\begin{document}
\maketitle

\begin{abstract}
Almost all materials are anisotropic. 
In this paper,  interface relations of anisotropic elliptic partial differential equations  involving discontinuities across  interfaces  are derived in two and three dimensions. Compared with  isotropic cases, the invariance of partial differential equations and the jump conditions under orthogonal coordinates transformation is not valid anymore. A systematic approach to derive the interface relations is established in this paper for anisotropic elliptic interface problems, which can be important for deriving high order accurate numerical methods.  

\end{abstract}

{\bf keywords:}
anisotropic elliptic interface problems, jump conditions/relations, local coordinates. 

 {\bf AMS Subject Classification 2000}
 .


\section{Introduction}

\smallskip

Consider an anisotropic elliptic interface problem below, 
\eqml{pde}
   && \dsp -\nabla \cdot \left (\mathbf{A} (\mathbf{x}) \nabla u ( \mathbf{x})  \right ) + \sigma (\mathbf{x}) u ( \mathbf{x}) =f(\mathbf{x}),   \qquad  \mathbf{x}\in \Omega\setminus \Gamma, \quad \Omega=\Omega^+ \cup \Omega^-, 
  \enml
where  $\mathbf{A}(\mathbf{x}) \in C(\Omega\setminus \Gamma)$ is a symmetric positive definite matrix whose eigenvalues satisfy, $\lambda_i  (\mathbf{x})\geq  \lambda_0>0$,  for all $\mathbf{x}\in \Omega$;  and $\sigma (\mathbf{x} )\geq 0$; the source term $f(\mathbf{x})\in C(\Omega\setminus \Gamma)$.  
%
\eqm
  [u] = w, \qquad [{\bf A}  \nabla u \cdot \mathbf{n} ]= v,
\enm
where $\Gamma\in C^2$ is a smooth interface with the solution domain, and the notation $[\;]$ is a jump of  a quantity across the interface. The above jump conditions are often referred as internal boundary conditions that make the problem well-posed. 

\begin{figure}[htbp]
\begin{minipage}[t]{2.5in}
$\null $ \includegraphics[width=1.05\textwidth]{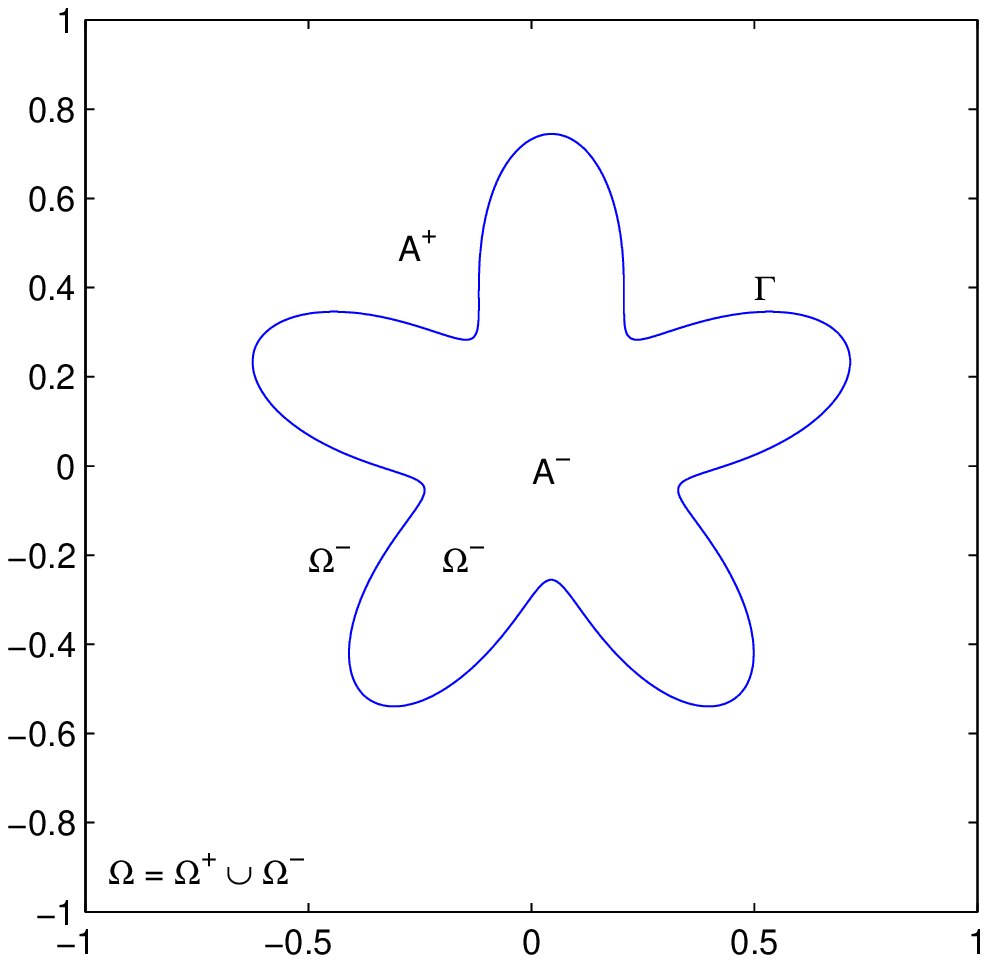}
\end{minipage}
\begin{minipage}[t]{2.9in}
\centering 
\includegraphics[width=0.65\textwidth]{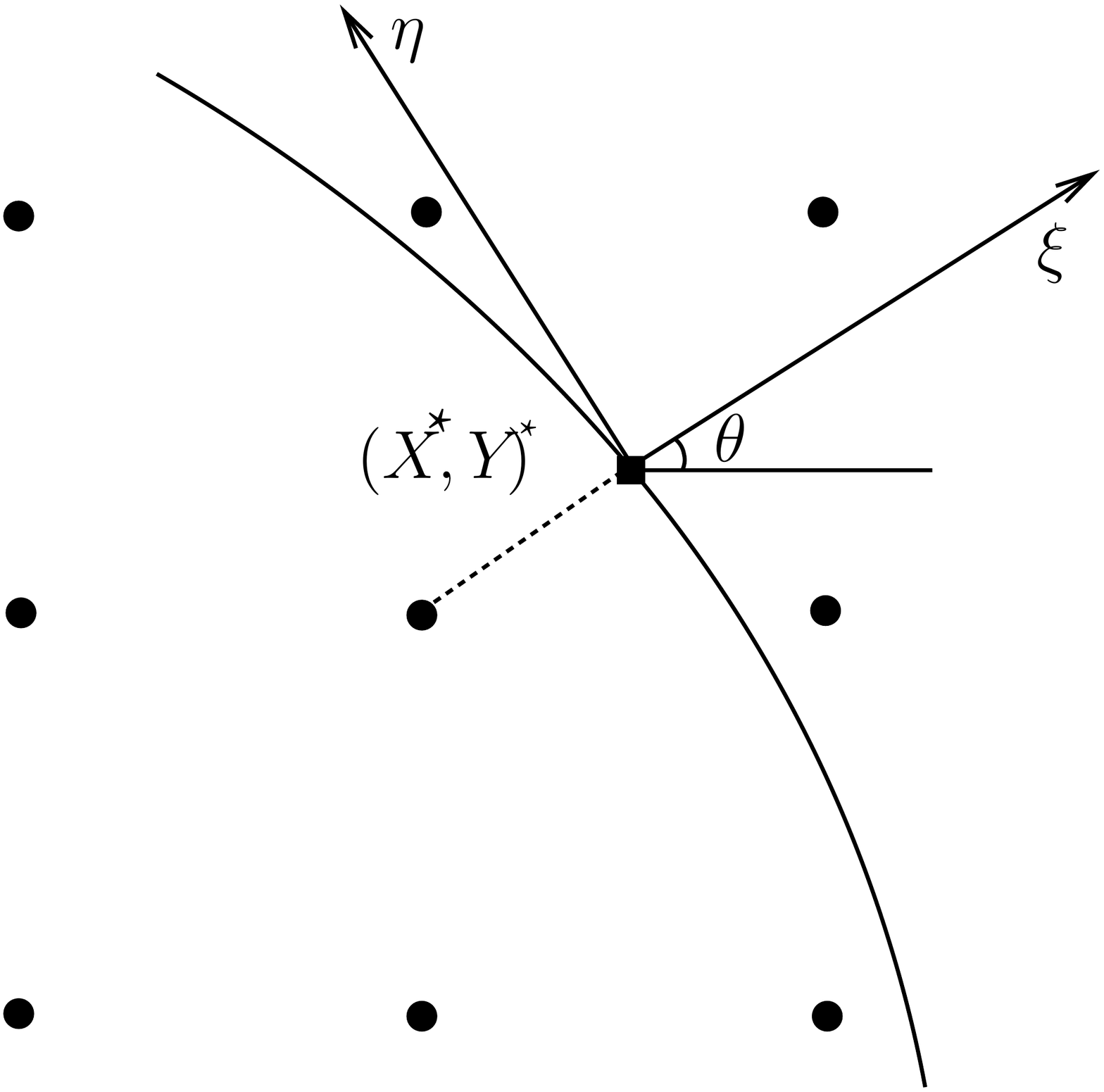}
 \end{minipage}

\caption{(a): A diagram of a domain $\Omega$ with an interface $\Gamma$. (b): A local coordinates.  }\label{domain}
\end{figure}

There are many discussions in the literature about anisotropic  interface problems including physics and modeling  
\cite{Bergmann17,Levitas-Warren16,McFadden93},  analysis  \cite{Suo90,Huang-Rokhlin92,Tuncel-Serbest},  and numerical methods  \cite{Dumett-Keener,An-Chen14,hou-wang-wang-sharp10}. 

For theoretical and numerical purposes, we want to know the interface relations of the partial derivatives of the solution across the interface.
When ${\bf A} = \beta {\bf I} $, that is, a scalar case,   the interface relations have been derived using the invariance of the PDE and the jump conditions. The derived interface relations have then been applied in derive accurate  finite difference methods such as the IIM \cite{rjl-li:sinum,deng-li,li:book}.

\section{The interface relations for anisotropic interface problems in 2D}

In 2D, the interface $\Gamma$ is a curve within the solution domain. We cannot assume invariances of the PDE, the flux jump conditions, and their surface derivatives under different coordinates systems. In this paper, we propose a way to parameterize the interface locally using  a level set function so that we can derive the interface relations in a systematically way as described below. 

\subsection{The local coordinate system and representation of the interface in 2D}

Let $(X^*, Y^*)$ be a fixed point on the interface $\Gamma$, and the normal direction at $(X^*, Y^*)$  be $(\cos \theta^*, \sin \theta^*)$, where $\theta^*$ is the angle of the normal direction and the $x$-axis, see Figure~\ref{domain} for an illustration.  The local coordinates in the neighborhood of $(X^*, Y^*)$  is defined as
\eqml{coor}
 \left \{ \begin{array}{l}
\xi =  \dsp ( x - X^* ) \, \mbox{cos}\, \theta^* + ( y - Y^* ) \, \mbox{sin}
\, \theta^* , \\ \eqsp
\eta  = -( x - X ^*) \,\mbox{sin}\,\theta ^* + ( y - Y^* ) \, \mbox{cos}\, \theta^*.
\end{array} \right.
\end{array}
\end{equation}
In a neighborhood of $(X^*,Y^*)$, the interface can be written as
\eqm
  \xi=\chi(\eta), \quad \mbox{with} \quad \chi(0)=0, \quad \chi'(0)=0,
\enm
and can be further parameterized as 
%
\eqml{coor2}
 \left \{ \begin{array}{l}
 X(s) =  \dsp   X^*  + \chi(s)  \, \mbox{cos}\, \theta^*  - s  \, \mbox{sin}
\, \theta^*, \\ \eqsp
Y(s)  =  Y^* +  \chi(s)  \,\mbox{sin}\,\theta ^* +  s  \, \mbox{cos}\, \theta^*.
\end{array} \right.
\end{array}
\end{equation}
We  have  $(X(0),Y(0))  = (X^*,   Y^*)$, here $s$ can be regarded as arc-length parameter starting from $(X^*, Y^*)$. The tangent vector then is
\eqml{tang1}
 \bftau (s)  &= & \dsp \left [  \frac{X'}{\sqrt{(X')^{2} + (Y')^{2} }}, \;  \frac{Y'}{\sqrt{(X')^{2} + (Y')^{2} }}\right ] \\ \eqsp
   &= & \dsp \left [  \frac{\chi' \, \cos \theta^* -\sin \theta^* }{\sqrt{ 1 +  (\chi ')^{2} }}, \;  \frac{ \chi' \, \sin \theta^* + \cos \theta^*}{\sqrt{ 1 +  (\chi ')^{2} }}\right ],
\enml
with  $\bftau (0)= [ - \sin \theta^*, \; \cos \theta^*]$,  and 
the normal direction is
\eqml{tang2}
\mathbf{n} = \dsp \left [  \frac{ \chi' \, \sin \theta^* + \cos \theta^*  }{\sqrt{ 1 +  (\chi ')^{2} }},  \;  \frac{ -\chi' \, \cos \theta^* + \sin \theta^* }{\sqrt{ 1 +  (\chi ')^{2} }}
\right ].
\enml

For simplicity, we still use the same  notations for the solution $u$, ${\bf A}$, $[u]=w(s)$, and $[{\bf A}\grad u \cdot {\bf n}] =v(s)$ in the local
coordinate system. In the neighborhood of $(X^*, Y^*)$, using the idea of the level set method, we can extend the quantities on the interface along the normal line using $\varphi(\xi, s)= \xi - \chi(s)$. Thus the tangential and normal derivatives and other interface quantities are  also defined in the neighborhood as the normal extension of  their value from the interface along the normal line.
Note that in the local coordinates, the PDE can be written as 
\begin{equation}\label{e:Formulas3-2}
a_{11}u_{\xi\xi}+2a_{12}u_{\xi\eta}+a_{22}u_{\eta\eta} - \sigma  u = f ,
\end{equation}
where $a_{ij}$ are defined below
\eqml{Aijtoaij}
a_{11} &=& \dsp A_{11}\cos ^{2} \theta^* + A_{22}\sin ^{2} \theta^* + 2A_{12}\cos \theta^* \sin\theta^*,\\ \eqsp
a_{22} &=& \dsp A_{11}\sin ^{2} \theta^* + A_{22}\cos ^{2} \theta^* - 2A_{12}\cos \theta^* \sin\theta^*,\\ \eqsp
a_{12} &=& \dsp  \left (A_{22}-A_{11} \right )\cos \theta^* \sin \theta^* + A_{12} \left (\cos ^{2}  \theta^*- \sin ^{2} \theta^* \right).
\enml
Under the framework above, we are ready to prove the main theorem in 2D.
\begin{theorem}
 If ${\bf A}$ is a piecewise constant matrix,  $u(x,y)\in C^2(\Omega^{\pm})$, $f(x,y)\in C(\Omega^{\pm})$, $\Gamma\in C^2$, $w\in C^2$, $v\in C^1$,  then the following interface relations hold.
\end{theorem}
\eqml{jump2}
&& u^{+}=u^{-}+w \\ \eqsp 
&& \dsp u_{\xi}^{+}=\frac{a^{-}_{11}}{a^{+}_{11}}
u_{\xi}^{-}-\frac{[a_{12}]}{a^{+}_{11}}u_{\eta}^{-}+\frac{1}{a^{+}_{11}}v- \frac{a^{+}_{12}}{a^{+}_{11}}w^{\prime} \\ \eqsp
&& \dsp u_{\eta}^{+}=u_{\eta}^{-}+w^{\prime} \\ \eqsp 
&& \dsp u_{\eta\eta}^{+}=\chi^{\prime\prime}\frac{[a_{11}]}{a^{+}_{11}}u_{\xi}^{-}+\chi^{\prime\prime} \frac{[a_{12}]}{a^{+}_{11}}u_{\eta}^{-}+u_{\eta\eta}^{-}- \frac{\chi^{\prime\prime}}{a^{+}_{11}}v+\chi^{\prime\prime}\frac{a^{+}_{12}}{a^{+}_{11}}w^{\prime}+w^{\prime \prime} \\ \eqsp
&&  \dsp u_{\xi \eta}^{+} = \chi^{\prime\prime}\frac{a^{+}_{11}[a_{12}]-2a^{+}_{12}[a_{11}]} {(a^{+}_{11})^{2}} u_{\xi}^{-}+ \chi^{\prime\prime}\frac{a^{+}_{11}[a_{22}]-2a^{+}_{12}[a_{12}]} {(a^{+}_{11})^{2}} u_{\eta}^{-} +\frac{a^{-}_{11}}{a^{+}_{11}}u_{\xi\eta}^{-} \\ \eqsp 
&& \dsp \null \quad - \frac{[a_{12}]}{a^{+}_{11}}u_{\eta\eta}^{-} +\chi^{\prime\prime} \frac{2a^{+}_{12}}{(a^{+}_{11})^2}v+\frac{1}{a^{+}_{11}}v^{\prime} +\chi^{\prime\prime}\frac{(a^{+}_{11}a^{+}_{22}-2(a^{+}_{12})^{2})} {(A^{+}_{11})^2}w^{\prime}  -\frac{a^{+}_{12}}{a^{+}_{11}}w^{\prime\prime} \\ \eqsp
 && \dsp  u_{\xi\xi}^{+} = - S_1 \chi^{\prime\prime}  u_{\xi}^{-} - S_2  \chi^{\prime\prime} u_{\eta}^{-} +\frac{a_{11}^{-}}{a_{11}^{+}} u_{\xi \xi}^{-}  
    +      \frac{2 a_{12}^{+}[a_{12}]-a_{11}^{+}[a_{22}]}{(a_{11}^{+})^{2}}u_{\xi \eta}^{-}  \\ \eqsp
    && \dsp  \null \qquad \qquad + \frac{2 (a_{12}^{+}[a_{11}]-a_{11}^{+}[a_{12}])} {(a_{11}^{+})^{2}}u_{\eta\eta}^{-} + \frac{[f]}{a_{11}^{+}} 
    + \frac {  [\sigma]} { a^{+}_{11} } u^{-}  + S_3,
\enml  
 where $S_1$, $S_2$, and   $S_3$ are given below,
\eqml {uintjB}
 S_1 &=&  \dsp  \frac{2a_{11}^{+}a_{12}^{+}[a_{12}]-4(a_{12}^{+})^{2}[a_{11}]+ a_{11}^{+}a_{22}^{+}[a_{11}]}{(a_{11}^{+})^{3}}  \\ \eqsp
 S_2  &=&  \dsp -\chi^{\prime\prime} \frac{2a_{11}^{+}a_{12}^{+}[a_{22}]-4(a_{12}^{+})^{2}[a_{12}]+ a_{11}^{+}a_{22}^{+}[a_{12}]}{(a_{11}^{+})^{3}} \\ \eqsp
 S_3  &=& \dsp  -\chi^{\prime\prime}\frac{4(a_{12}^{+})^{2}-a_{11}^{+}a_{22}^{+}}{(a_{11}^{+})^{3}}v
-\frac{2a_{12}^{+}}{(a_{11}^{+})^{2}}v^{\prime} \\ \eqsp
 &&  \null \quad   
 \dsp + \frac { \sigma^{+} w} { a^{+}_{11} }   - \chi^{\prime\prime}\frac{a_{11}^{+}a_{12}^{+}a_{22}^{+}-4(a_{12}^{+})^{3}} {(a_{11}^{+})^{3}} w^{\prime} + \frac{2(a_{12}^{+})^{2}-a_{11}^{+}a_{12}^{+}}{(a_{11}^{+})^{2}} w^{\prime\prime} .
\enml
Note that $S_3$ depends on the jump conditions $w$ and $v$, their (surface) derivatives,  coefficient matrix ${\bf A}$, and the curvature of $\Gamma$. 

\ \

{\bf Proof:} Differentiating $[u]=w$  with respect to $s$  once we get
\eqm
[u_{\xi}]\;\chi' + [u_{\eta}] = w'(s).
\enm
%
%
Differentiating  the identity above with respect to $s$  we get
\begin{equation}\label{e:Formulas3-11}
[u_{\xi \xi}]\chi^{\prime 2}+2[u_{\xi \eta}]\chi^{\prime}+[u_{\xi}]\chi^{\prime \prime}+[u_{\eta \eta}]=w^{\prime \prime}(s).
\end{equation}
Setting $s=0$ in the above two identities and using $\chi^{\prime}(0)=0$, we obtain the third and the  fourth identities in the theorem. 

From $u_{x}=u_{\xi}\cos \theta^* - u_{\eta}\sin \theta^*$, $u_{y}=u_{\xi}\sin \theta^* + u_{\eta}\cos \theta^*$ and $\mathbf{n(s)}$ defined above, and the relation between $A_{ij}$ and $a_{ij}$, and with some manipulations, 
we can rewrite $ \big[ \mathbf{A} \nabla u \cdot \mathbf{n} \big]$~as,
\[
\begin{aligned}
& \big[ \mathbf{A} \nabla u \cdot \mathbf{n} \big]
=
\left[
\left(
      \begin{array}{cc}
        A_{11} & A_{12} \\
        A_{12} & A_{22} \\
      \end{array}
\right)
\left(
      \begin{array}{c}
      u_{x} \\
      u_{y} \\
      \end{array}
\right) \cdot
\left(
      \begin{array}{c}
      n_{x} \\
      n_{y} \\
      \end{array}
\right)
\right]\\
& = 
\bigg[\left(
      \begin{array}{c}
      A_{11}(u_{\xi}\cos \theta^* - u_{\eta}\sin \theta^*) + A_{12}(u_{\xi}\sin \theta^* + u_{\eta}\cos \theta^*)\\
      A_{12}(u_{\xi}\cos \theta^* - u_{\eta}\sin \theta^*) + A_{22}(u_{\xi}\sin \theta^* + u_{\eta}\cos \theta^*)\\
      \end{array}
\right)^{T} \cdot \mathbf{n}  
\bigg]\\ \eqsp
&= \frac{1}{\sqrt{1+(\chi^{ \prime})^{2}}} \left [ \frac{\null}{\null} \left (a_{11}-\chi^{ \prime} a_{12} \right )u_{\xi}] +  \left [(a_{12}-\chi^{ \prime} a_{22}
 \right ) u_{\eta} \right].
\end{aligned}
\]
Therefore, the flux jump condition can be  written as
\begin{equation}\label{e:Formulas3-14}
 \left [  \frac{\null}{\null} \left  (a_{11}-\chi^{ \prime} a_{12} \right )u_{\xi} \right ] + \left [ \frac{\null}{\null} \left (a_{12}-\chi^{ \prime} a_{22} \right )u_{\eta} \right ]
 =\sqrt{1+(\chi^{ \prime})^{2}} \,v(s).
\end{equation}
We get the second identity by setting $s=0$. 

To get the fifth identity, we  differentiate the flux jump condition above  with respect to~$s$. The left hand side then is
\eqmno
  LHS &=&  \left [  \frac{\null}{\null} \left  (a_{11}-\chi^{ \prime} a_{12} \right )  \left ( u_{\xi \xi} \chi' + u_{\xi \eta}  \right )  \right ]  +
   \left [  \frac{\null}{\null}  - \chi''  a_{12} u_{\xi} \right ]  \\ \eqsp
   &&  \null +  \left [ \frac{\null}{\null} \left (a_{12}-\chi^{ \prime} a_{22} \right ) \left ( u_{\eta\xi} \chi' + u_{\eta\eta}  \right )\right ]
 +  \left [ \frac{\null}{\null} -\chi^{ \prime  \prime} a_{22}  u_{\eta} \right ]. 
\enmno
The right hand side is
\eqmno
 RHS = \sqrt{1+(\chi^{ \prime})^{2}} \,v^{\prime} (s) + \frac{ \chi^{ \prime} \chi^{ \prime \prime}}{\sqrt{1+(\chi^{ \prime})^{2}}} \, v(s).
\enmno
By plugging $s=0$ to the left and right hand sides and using $\chi^{ \prime} (0)=0$, we get the fifth identity. 

The last identity is obtained from the PDE using
\begin{equation}\label{e:Formulas3-17}
[a_{11}u_{\xi\xi}]+2[a_{12}u_{\xi\eta}]+[a_{22}u_{\eta\eta}] - [\sigma u]=[f].
\end{equation}

\begin{remark}
 For variable coefficients ${\bf A({\bf x}})$,  $\sigma({\bf x})$, and $f({\bf x})$, the derivation process is similar but the expressions are long and more complicated.  The PDE in the local coordinates now is
 \eqm
   a_{11}u_{\xi\xi} + 2 a_{12} u_{\xi\eta}  + a_{22} u_{\eta\eta} +  c_{1} u_{\xi}  + c_{2}u_{\eta} - \sigma u = f, 
 \enm
 where
 \eqm
   c_{1}=  \frac{ \partial a_{11}}{\partial \xi}  +  \frac{ \partial a_{12}}{\partial \eta} , \qquad   c_{2}  =  \frac{ \partial a_{12}}{\partial \xi}  +  \frac{ \partial a_{22}}{\partial \eta} 
 \enm
 The first four  identities are the same and the last two  are the following: 
\end{remark}
\eqml{jumpvar}
 u^{+}_{\xi\eta}&=& \dsp \frac{c^{-}_{3}[a_{11}]-a_{11}^{-}[c_{3}]-\chi^{\prime \prime}a^{+}_{12}[a_{11}]}{(a_{11}^{+})^2}u^{-}_{\xi} + \frac{c^{+}_{3}[a_{12}]-a_{11}^{+}[c_{4}]-\chi^{\prime \prime}a^{+}_{12}[a_{12}]}{(a_{11}^{+})^2}u^{-}_{\eta} - \frac{[a_{12}]}{a_{11}^{+}}u^{-}_{\eta\eta}  \\ \eqsp
&& \dsp \null  + \frac{a_{11}^{-}}{a_{11}^{+}}u^{-}_{\xi\eta} + \frac{c^{+}_{3}a^{+}_{12}-a_{11}^{+}c^{+}_{4}-\chi^{\prime \prime}(a^{+}_{12})^{2}}{(a_{11}^{+})^2}w^{\prime} - \frac{a_{12}^{+}}{a_{11}^{+}}w^{\prime\prime} 
+ \frac{\chi^{\prime \prime}a^{+}_{12}-c^{+}_{3}}{(a_{11}^{+})^2}v + \frac{v^{\prime}}{a_{11}^{+}},
\enml
\eqml{jumpvar}
u^{+}_{\xi\xi} &=& \dsp 
 \frac{(a_{11}^{+})^{2}c^{-}_{1} - a_{11}^{+}a_{11}^{-}c^{+}_{1} - 2a^{+}_{12}(c^{-}_{3}[a_{11}]-a_{11}^{-}[c_{3}]-\chi^{\prime \prime}a^{+}_{12}[a_{11}])- \chi^{\prime \prime}a^{+}_{11}a^{+}_{22}[a_{11}]}{(a_{11}^{+})^3}u^{-}_{\xi}\\ \eqsp
&& \dsp  \null  + \frac{ a_{11}^{+}c^{+}_{1}[a_{12}] -(a_{11}^{+})^{2}[c_{2}] - 2a^{+}_{12}(c^{+}_{3}[a_{12}]-a_{11}^{+}[c_{4}]-\chi^{\prime \prime}a^{+}_{12}[a_{12}])- \chi^{\prime \prime}a^{+}_{11}a^{+}_{22}[a_{12}]}{(a_{11}^{+})^3}u^{-}_{\eta}\\ \eqsp
&& \dsp  \null  + \frac{a_{11}^{-}}{a_{11}^{+}}u^{-}_{\xi\xi} + \frac{2a^{+}_{12}[a_{12}]-a^{+}_{11}[a_{22}]}{(a_{11}^{+})^2}u^{-}_{\eta\eta} + \frac{2(a^{+}_{12}[a_{11}]-a^{+}_{11}[a_{12}])}{(a_{11}^{+})^2}u^{-}_{\xi\eta}\\ \eqsp
&& \dsp \null  + \frac{a_{11}^{+}a^{+}_{12}c^{+}_{1} - (a_{11}^{+})^{2}c^{+}_{2} -2a^{+}_{12}(c^{+}_{3}a^{+}_{12}-a_{11}^{+}c^{+}_{4}-\chi^{\prime \prime}(a^{+}_{12})^{2}) -\chi^{\prime \prime} a_{11}^{+}a_{12}^{+}a_{22}^{+}}{(a_{11}^{+})^3}w^{\prime}\\ \eqsp
&& \dsp \null   + \frac{2(a_{12}^{+})^{2}-a_{11}^{+}a_{22}^{+}}{a_{11}^{+}}w^{\prime\prime} + \frac{\chi^{\prime \prime}a^{+}_{11}a^{+}_{22}-a^{+}_{11}c^{+}_{1} - 2a^{+}_{12}(\chi^{\prime \prime}a^{+}_{12}- c^{+}_{3})}{(a_{11}^{+})^3}v - \frac{2a^{+}_{12}}{(a_{11}^{+})^2}v^{\prime}\\ \eqsp
&& \dsp  \null  + \frac{\sigma^{+}[u]+[\sigma]u^{-}}{a_{11}^{+}} + \frac{[f]}{a_{11}^{+}},
\enml
where $c_k$, $k=1,2,3,4$ are given below
where
 \eqm
   c_{1}=  \frac{ \partial a_{11}}{\partial \xi}  +  \frac{ \partial a_{12}}{\partial \eta} , \quad   c_{2}  =  \frac{ \partial a_{12}}{\partial \xi}  +  \frac{ \partial a_{22}}{\partial \eta} , \quad c_{3}=\frac{ \partial a_{11}}{\partial \eta}  -\chi^{\prime \prime}a_{12}, \quad c_{4}=\frac { \partial a_{12} } { \partial \eta } - \chi^{ \prime \prime } a_{22}.
 \enm

\parskip= 0.0cm
\renewcommand{\baselinestretch}{1.0}

\bibliographystyle{amsplain}

\bibliography{../../BIB/bib,../../BIB/zhilin,anis,../../BIB/other}
\end{document}